\documentclass[12pt]{article}

\usepackage{amsmath,amssymb,amsthm,amscd,a4wide}
\usepackage{hyperref}
\hypersetup{colorlinks,citecolor=blue,filecolor=black,linkcolor=blue,urlcolor=blue}

\begin{document}

\newcommand{\ad}{{\rm ad}}
\newcommand{\cri}{{\rm cri}}
\newcommand{\End}{{\rm{End}\ts}}
\newcommand{\Rep}{{\rm{Rep}\ts}}
\newcommand{\Hom}{{\rm{Hom}}}
\newcommand{\Mat}{{\rm{Mat}}}
\newcommand{\ch}{{\rm{ch}\ts}}
\newcommand{\chara}{{\rm{char}\ts}}
\newcommand{\diag}{{\rm diag}}
\newcommand{\non}{\nonumber}
\newcommand{\wt}{\widetilde}
\newcommand{\wh}{\widehat}
\newcommand{\ot}{\otimes}
\newcommand{\la}{\lambda}
\newcommand{\La}{\Lambda}
\newcommand{\De}{\Delta}
\newcommand{\al}{\alpha}
\newcommand{\be}{\beta}
\newcommand{\ga}{\gamma}
\newcommand{\Ga}{\Gamma}
\newcommand{\ep}{\epsilon}
\newcommand{\ka}{\kappa}
\newcommand{\vk}{\varkappa}
\newcommand{\si}{\sigma}
\newcommand{\vs}{\varsigma}
\newcommand{\vp}{\varphi}
\newcommand{\de}{\delta}
\newcommand{\ze}{\zeta}
\newcommand{\om}{\omega}
\newcommand{\Om}{\Omega}
\newcommand{\ee}{\epsilon^{}}
\newcommand{\su}{s^{}}
\newcommand{\hra}{\hookrightarrow}
\newcommand{\ve}{\varepsilon}
\newcommand{\ts}{\,}
\newcommand{\vac}{\mathbf{1}}
\newcommand{\di}{\partial}
\newcommand{\qin}{q^{-1}}
\newcommand{\tss}{\hspace{1pt}}
\newcommand{\Sr}{ {\rm S}}
\newcommand{\U}{ {\rm U}}
\newcommand{\BL}{ {\overline L}}
\newcommand{\BE}{ {\overline E}}
\newcommand{\BP}{ {\overline P}}
\newcommand{\AAb}{\mathbb{A}\tss}
\newcommand{\CC}{\mathbb{C}\tss}
\newcommand{\KK}{\mathbb{K}\tss}
\newcommand{\QQ}{\mathbb{Q}\tss}
\newcommand{\SSb}{\mathbb{S}\tss}
\newcommand{\ZZ}{\mathbb{Z}\tss}
\newcommand{\X}{ {\rm X}}
\newcommand{\Y}{ {\rm Y}}
\newcommand{\Z}{{\rm Z}}
\newcommand{\Ac}{\mathcal{A}}
\newcommand{\Lc}{\mathcal{L}}
\newcommand{\Mc}{\mathcal{M}}
\newcommand{\Pc}{\mathcal{P}}
\newcommand{\Qc}{\mathcal{Q}}
\newcommand{\Rc}{\mathcal{R}}
\newcommand{\Sc}{\mathcal{S}}
\newcommand{\Tc}{\mathcal{T}}
\newcommand{\Bc}{\mathcal{B}}
\newcommand{\Ec}{\mathcal{E}}
\newcommand{\Fc}{\mathcal{F}}
\newcommand{\Gc}{\mathcal{G}}
\newcommand{\Hc}{\mathcal{H}}
\newcommand{\Uc}{\mathcal{U}}
\newcommand{\Vc}{\mathcal{V}}
\newcommand{\Wc}{\mathcal{W}}
\newcommand{\Yc}{\mathcal{Y}}
\newcommand{\Ar}{{\rm A}}
\newcommand{\Br}{{\rm B}}
\newcommand{\Ir}{{\rm I}}
\newcommand{\Fr}{{\rm F}}
\newcommand{\Jr}{{\rm J}}
\newcommand{\Or}{{\rm O}}
\newcommand{\GL}{{\rm GL}}
\newcommand{\Spr}{{\rm Sp}}
\newcommand{\Rr}{{\rm R}}
\newcommand{\Zr}{{\rm Z}}
\newcommand{\gl}{\mathfrak{gl}}
\newcommand{\middd}{{\rm mid}}
\newcommand{\ev}{{\rm ev}}
\newcommand{\Pf}{{\rm Pf}}
\newcommand{\Norm}{{\rm Norm\tss}}
\newcommand{\oa}{\mathfrak{o}}
\newcommand{\spa}{\mathfrak{sp}}
\newcommand{\osp}{\mathfrak{osp}}
\newcommand{\f}{\mathfrak{f}}
\newcommand{\g}{\mathfrak{g}}
\newcommand{\h}{\mathfrak h}
\newcommand{\n}{\mathfrak n}
\newcommand{\z}{\mathfrak{z}}
\newcommand{\Zgot}{\mathfrak{Z}}
\newcommand{\p}{\mathfrak{p}}
\newcommand{\sll}{\mathfrak{sl}}
\newcommand{\agot}{\mathfrak{a}}
\newcommand{\qdet}{ {\rm qdet}\ts}
\newcommand{\Ber}{ {\rm Ber}\ts}
\newcommand{\HC}{ {\mathcal HC}}
\newcommand{\cdet}{{\rm cdet}}
\newcommand{\rdet}{{\rm rdet}}
\newcommand{\tr}{ {\rm tr}}
\newcommand{\gr}{{\rm gr}}
\newcommand{\str}{ {\rm str}}
\newcommand{\loc}{{\rm loc}}
\newcommand{\Gr}{{\rm G}}
\newcommand{\sgn}{ {\rm sgn}\ts}
\newcommand{\sign}{{\rm sgn}}
\newcommand{\ba}{\bar{a}}
\newcommand{\bb}{\bar{b}}
\newcommand{\bi}{\bar{\imath}}
\newcommand{\bj}{\bar{\jmath}}
\newcommand{\bk}{\bar{k}}
\newcommand{\bl}{\bar{l}}
\newcommand{\hb}{\mathbf{h}}
\newcommand{\Sym}{\mathfrak S}
\newcommand{\fand}{\quad\text{and}\quad}
\newcommand{\Fand}{\qquad\text{and}\qquad}
\newcommand{\For}{\qquad\text{or}\qquad}
\newcommand{\OR}{\qquad\text{or}\qquad}

\renewcommand{\theequation}{\arabic{section}.\arabic{equation}}

\newtheorem{thm}{Theorem}[section]
\newtheorem{lem}[thm]{Lemma}
\newtheorem{prop}[thm]{Proposition}
\newtheorem{cor}[thm]{Corollary}
\newtheorem{conj}[thm]{Conjecture}
\newtheorem*{mthm}{Main Theorem}
\newtheorem*{mthma}{Theorem A}
\newtheorem*{mthmb}{Theorem B}
\newtheorem*{mthmc}{Theorem C}
\newtheorem*{mthmd}{Theorem D}
\newtheorem*{ffthm}{Feigin--Frenkel theorem}

\theoremstyle{definition}
\newtheorem{defin}[thm]{Definition}

\theoremstyle{remark}
\newtheorem{remark}[thm]{Remark}
\newtheorem{example}[thm]{Example}

\newcommand{\bth}{\begin{thm}}
\renewcommand{\eth}{\end{thm}}
\newcommand{\bpr}{\begin{prop}}
\newcommand{\epr}{\end{prop}}
\newcommand{\ble}{\begin{lem}}
\newcommand{\ele}{\end{lem}}
\newcommand{\bco}{\begin{cor}}
\newcommand{\eco}{\end{cor}}
\newcommand{\bde}{\begin{defin}}
\newcommand{\ede}{\end{defin}}
\newcommand{\bex}{\begin{example}}
\newcommand{\eex}{\end{example}}
\newcommand{\bre}{\begin{remark}}
\newcommand{\ere}{\end{remark}}
\newcommand{\bcj}{\begin{conj}}
\newcommand{\ecj}{\end{conj}}

\newcommand{\bal}{\begin{aligned}}
\newcommand{\eal}{\end{aligned}}
\newcommand{\beq}{\begin{equation}}
\newcommand{\eeq}{\end{equation}}
\newcommand{\ben}{\begin{equation*}}
\newcommand{\een}{\end{equation*}}

\newcommand{\bpf}{\begin{proof}}
\newcommand{\epf}{\end{proof}}

\def\beql#1{\begin{equation}\label{#1}}

\title{\Large\bf Segal--Sugawara vectors for the Lie algebra of type $G_2$}

\author{{A. I. Molev,\quad E. Ragoucy\quad and\quad N. Rozhkovskaya}}

\date{} 
\maketitle

\vspace{20 mm}

\begin{abstract}
Explicit formulas for Segal--Sugawara vectors associated with
the simple Lie algebra $\g$ of type $G_2$ are found by using computer-assisted
calculations. This leads to a direct proof of the
Feigin--Frenkel theorem describing the center of the
corresponding affine vertex algebra at the critical level.
As an application, we give an explicit solution of Vinberg's quantization
problem by providing formulas for generators of maximal
commutative subalgebras of $\U(\g)$. We also calculate the eigenvalues
of the Hamiltonians on the Bethe vectors in the Gaudin model associated with $\g$.

\medskip

Preprint LAPTH-006/16

\end{abstract}


\vspace{30 mm}

\noindent
School of Mathematics and Statistics\newline
University of Sydney,
NSW 2006, Australia\newline
alexander.molev@sydney.edu.au

\vspace{7 mm}

\noindent
Laboratoire de Physique Th\'{e}orique LAPTh,
CNRS and Universit\'{e} de Savoie\newline
BP 110, 74941 Annecy-le-Vieux Cedex, France\newline
eric.ragoucy@lapth.cnrs.fr

\vspace{7 mm}

\noindent
Department of Mathematics\newline
Kansas State University, USA\newline
rozhkovs@math.ksu.edu

\newpage

%

\section{Introduction}
\label{sec:int}
\setcounter{equation}{0}

For a simple Lie algebra $\g$ over $\CC$
equipped with
a standard symmetric invariant bilinear form,
consider the corresponding (non-twisted)
affine Kac--Moody algebra $\wh\g$
\beql{km}
\wh\g=\g\tss[t,t^{-1}]\oplus\CC K.
\eeq
The universal vacuum module $V(\g)$ over $\wh\g$
is the quotient of $\U(\wh\g)$ by the left ideal generated by $\g[t]$.
A {\em Segal--Sugawara vector} is any element $S\in V(\g)$ with the property
\ben
\g[t]\tss S\in (K+h^{\vee})\ts V(\g),
\een
where $h^{\vee}$ is the dual Coxeter number for $\g$. In particular, the canonical quadratic
Segal--Sugawara vector is given by
\beql{canss}
S=\sum_{a=1}^d X_a[-1]^2,
\eeq
where $X_1,\dots,X_d$ is an orthonormal basis of
$\g$ and we write $X[r]=X\tss t^r$ for $X\in\g$.

By an equivalent approach, Segal--Sugawara vectors are elements of the subspace
$\z(\wh\g)$ of invariants of the vacuum module at the critical level
\beql{ffvac}
\z(\wh\g)=\{v\in V(\g)_{\text{cri}}\ |\ \g[t]\tss v=0\},
\eeq
where $V(\g)_{\text{cri}}$ is the quotient of $V(\g)$ by the submodule
$(K+h^{\vee})\ts V(\g)$.
The vacuum module possesses a vertex algebra structure,
and by the definition \eqref{ffvac},
$\z(\wh\g)$ coincides with the {\it center\/}
of the vertex algebra $V(\g)_{\text{cri}}$.
This induces a structure of commutative associative algebra on the center
which coincides with the one obtained via identification of $\z(\wh\g)$
with a subalgebra of $\U\big(t^{-1}\g[t^{-1}]\big)$.

The structure of $\z(\wh\g)$ was
described by a theorem of Feigin and Frenkel~\cite{ff:ak}, which states that
$\z(\wh\g)$ is an algebra of polynomials in infinitely many variables; see \cite{f:lc} for
a detailed exposition of these results. The algebra $\z(\wh\g)$ is refereed to as
the {\em Feigin--Frenkel center}.
Explicit formulas for generators
of this algebra were found in \cite{ct:qs} for type $A$ and in \cite{m:ff}
for types $B$, $C$ and $D$; see also \cite{cm:ho} and \cite{mr:mm} for simpler arguments
in type $A$ and extensions to Lie superalgebras.

Our goal in this paper is to give explicit formulas for generators of $\z(\wh\g)$
in the case where $\g$ is the exceptional Lie algebra of type $G_2$.
In particular, we obtain a direct proof of the Feigin--Frenkel theorem in this case.
Furthermore, using the connections with the Gaudin model as discovered in \cite{ffr:gm},
we get formulas for higher Gaudin Hamiltonians associated with $\g$ and calculate
their eigenvalues on the Bethe vectors; see also \cite{fft:gm}.
In the classical types such formulas were given in a recent work \cite{mm:eb}.

As another application, following \cite{fft:gm} and
\cite{r:si} we give explicit formulas for algebraically independent generators of
maximal commutative subalgebras of $\U(\g)$. These subalgebras $\Ac_{\mu}$
are parameterized by regular elements $\mu\in\g^*$, and their classical limits
$\overline\Ac_{\mu}$ are Poisson commutative subalgebras of $\Sr(\g)$ known as
the Mishchenko--Fomenko or shift of argument subalgebras.
The formulas for generators of $\Ac_{\mu}$ thus provide
an explicit solution of Vinberg's quantization
problem \cite{v:sc}.

Our calculations of the explicit expressions for the Segal--Sugawara vectors
and their Harish-Chandra images were computer-assisted. We gratefully acknowledge the use of the
{\em Symbolic Manipulation System FORM} originally developed by Vermaseren~\cite{v:nf}.

This project was completed
within the visitor program of the Center of Quantum Algebra
of the South China University of Technology, Guangzhou, China.
We would like to thank the Center and the School of Mathematical Sciences
for the warm hospitality during our visits.

\section{Lie algebra of type $G_2$ and its matrix presentation}
\label{sec:liea}
\setcounter{equation}{0}

We start by recalling well-known
matrix presentations of an arbitrary simple Lie algebra $\g$; see, e.g.,
\cite{d:qg} and \cite{g:ci}.

\subsection{Matrix presentations of simple Lie algebras}
\label{subsec:mps}

Equip $\g$ with
a symmetric invariant bilinear form $\langle\ts\ts,\ts\rangle$.
Choose a basis $X^1,\dots,X^d$ of $\g$ and let $X_1,\dots,X_d$ be its dual
with respect
to the form. Let $\pi$
be a faithful representation of $\g$ afforded by a finite-dimensional vector space $V$,
\beql{pi}
\pi:\g\to \End V.
\eeq
Introduce the elements
\beql{G}
G=\sum_{i=1}^d \pi(X^i)\ot X_i\in \End V\ot\U(\g)
\eeq
and
\beql{Om}
\Om=\sum_{i=1}^d \pi(X^i)\ot \pi(X_i)\in \End V\ot \End V.
\eeq
Note that $G$ and $\Om$ are independent of the choice of the basis $X^i$. In particular,
\beql{Omalt}
\Om=\sum_{i=1}^d \pi(X_i)\ot \pi(X^i).
\eeq
Consider the tensor product algebra $\End V\ot \End V\ot \U(\g)$
and identify $\Om$ with the element $\Om\ot 1$. Also, introduce its elements
\beql{Got}
G_1=\sum_{i=1}^d \pi(X^i)\ot 1\ot X_i\Fand
G_2=\sum_{i=1}^d 1\ot \pi(X^i)\ot X_i.
\eeq
Write the commutation relations for $\g$,
\beql{commg}
[X_i,X_j]=\sum_{k=1}^d c_{ij}^{\tss k}\ts  X_k
\eeq
with structure coefficients $c_{ij}^{\tss k}$. We will regard
the universal enveloping algebra $\U(\g)$ as the associative algebra
with generators $X_i$ subject to the defining relations \eqref{commg},
where the left hand side is understood as the commutator $X_iX_j-X_jX_i$.

\bpr\label{prop:matrrel}
The defining relations of $\U(\g)$ are equivalent to the matrix relation
\beql{matrelg}
G_1\tss G_2-G_2\tss G_1=-\Om\ts  G_2+G_2\tss \Om.
\eeq
\epr

\bpf
The left hand side of \eqref{matrelg} reads
\ben
\sum_{i,j=1}^d \pi(X^i)\ot \pi(X^j)\ot (X_iX_j-X_jX_i).
\een
For the right hand side we have
\beql{rhgom}
-\sum_{i,k=1}^d \pi(X^i)\ot \pi\big([X_i,X^k]\big)\ot X_k.
\eeq
By the invariance of the form, we find
\ben
\langle [X_i,X^k],X_j\rangle=-\langle X^k,[X_i,X_j]\rangle=-c_{ij}^{\tss k}.
\een
Hence \eqref{rhgom} equals
\ben
\sum_{i,j,k=1}^d c_{ij}^{\tss k}\ts\pi(X^i)\ot \pi(X^j)\ot X_k.
\een
Since the representation $\pi$ is faithful, we conclude
that \eqref{matrelg} is equivalent to
the defining relations \eqref{commg} of $\U(\g)$.
\epf

The defining relations \eqref{matrelg} can be written in an equivalent form
\beql{matrelgequi}
G_1\tss G_2-G_2\tss G_1=\Om\ts  G_1-G_1\tss \Om,
\eeq
which is easily verified with the use of \eqref{Omalt}.
The element $G$ can be regarded as an $n\times n$ matrix ($n=\dim V$) with
entries in $\U(\g)$. We point out another well-known property of this matrix
which goes back to \cite{g:tc}.

\bco\label{cor:casim}
All elements $\tr\ts G^k$ with $k\geqslant 1$ belong to the center of $\U(\g)$.
\eco

\bpf
Relation \eqref{matrelg} implies
\ben
G_1\tss G^k_2-G^k_2\tss G_1=-\Om\ts  G^k_2+G^k_2\tss \Om.
\een
By taking trace over the second copy of $\End V$ and using
its cyclic property, we get $[G_1,\tss \tr\ts G^k_2]=0$ as required.
\epf

The Casimir elements $\tr\ts G^k$ are widely used in representation theory,
especially for the Lie algebras $\g$ of classical types. In those cases
one usually takes $V$ to be the first fundamental (or vector) representation.

Note also that the relation \eqref{matrelg} can be regarded as the `classical part'
of the $RTT$ presentation of the Yangian $\Y(\g)$ associated with $\g$; see \cite{d:qg}.
More precisely, $\Y(\g)$ contains $\U(\g)$ as a subalgebra, and \eqref{matrelg}
is recovered as a reduction of the $RTT$ relation to the generators of this subalgebra.

\subsection{Lie algebra of type $G_2$}
\label{subsec:liegtwo}

The simple Lie algebra $\g$ of type $G_2$ admits a few different presentations;
see, e.g., \cite{fh:rt}, \cite{w:ccg}. It is well-known that it
can be embedded into the orthogonal Lie algebras $\oa_7$ and $\oa_8$; these embeddings
were employed in \cite{mr:cw} to construct the classical $\Wc$-algebra for $\g$.
We will follow \cite[Lect.~22]{fh:rt} to realize $\g$ as the direct sum
of vector spaces
\ben
\g=\CC^3\oplus\sll_3\oplus (\CC^3)^*.
\een
The Lie bracket on $\g$ is determined by the conditions that
$\sll_3$ is a subalgebra of $\g$, the vector spaces
$\CC^3$ and $(\CC^3)^*$ are, respectively, the vector representation
of $\sll_3$ and its dual, together with additional brackets
\ben
\CC^3\times \CC^3\to (\CC^3)^*,\qquad (\CC^3)^*\times (\CC^3)^*\to \CC^3
\Fand \CC^3\times (\CC^3)^*\to \sll_3.
\een
To produce a matrix presentation of $\g$ as provided by Proposition~\ref{prop:matrrel},
consider the $7$-dimensional representation $\pi:\g\to\End V$ with
$V\cong \CC^7$
where the action is described explicitly as follows. Write an arbitrary element
of $\g$ as a triple $(v,A,\vp)$, where $A\in \sll_3$ is a traceless $3\times 3$ matrix,
\ben
v=\begin{bmatrix}v_1\\v_2\\v_3\end{bmatrix}
\in \CC^3\Fand \vp=\big[\vp_1,\vp_2,\vp_3\big]
\in (\CC^3)^*.
\een
Then, under the representation $\pi$ we have
\ben
\pi:(v,A,\vp)\mapsto
\begin{bmatrix}A&v&\frac1{\sqrt2}B(\vp^t)\\[0.5em]
\vp&0&-v^t\\[0.5em]
\frac1{\sqrt2}B(v)&-\vp^t&-A^t
\end{bmatrix},
\een
where $t$ denotes the antidiagonal matrix transposition,
\ben
v^t=\big[v_3,v_2,v_1\big],\qquad
\vp^t=\begin{bmatrix}\vp_3\\\vp_2\\\vp_1\end{bmatrix},\qquad (A^t)_{ij}=A_{\tss 4-j,\tss 4-i}
\een
and
\ben
B(v)=\begin{bmatrix}\phantom{-}v_2&-v_1&\phantom{-}0\ts\\
-v_3&\phantom{-}0&\phantom{-}v_1\ts\\
\phantom{-}0&\phantom{-}v_3&-v_2\ts
\end{bmatrix}.
\een
The representation $\pi$ is faithful, so we may use it to identify
the Lie algebra $\g$ with its image under $\pi$ where elements of $\g$
can be regarded as $7\times 7$ matrices. Then the bilinear form on $\g$
defined by
\beql{formg}
\langle X,Y\rangle=\frac16\ts \tr\ts XY
\eeq
is symmetric and invariant. Note that this form is
proportional to the standard normalized Killing form
\ben
\frac{1}{2\tss h^{\vee}}\ts\tr\ts\big(\ad\tss X\ts\ad\tss Y\big),
\een
where $h^{\vee}=4$ is the dual Coxeter number for $\g$.
We have
\ben
\langle X,Y\rangle=\frac{1}{24}\ts\tr\ts\big(\ad\tss X\ts\ad\tss Y\big).
\een
The additional scalar factor $1/3$ is meant to simplify our formulas for
Segal--Sugawara vectors by avoiding fractions.

Let $e_{ij}\in\End\CC^7$ denote the standard
matrix units. For all $1\leqslant i,j\leqslant 7$
set $f_{ij}=e_{ij}-e_{j'i'}$, where $i'=8-i$. The following elements form
a basis of $\g$:
\ben
f_{11}-f_{22},\qquad f_{22}-f_{33},\qquad
f_{ij} \quad\text{with}\quad 1\leqslant i,j\leqslant 3\fand i\ne j,
\een
together with
\ben
f_{14}-\frac1{\sqrt2}\ts f_{3'2},\qquad
f_{24}-\frac1{\sqrt2}\ts f_{1'3},\qquad
f_{34}-\frac1{\sqrt2}\ts f_{2'1},
\een
and
\ben
f_{41}-\frac1{\sqrt2}\ts f_{23'},\qquad
f_{42}-\frac1{\sqrt2}\ts f_{31'},\qquad
f_{43}-\frac1{\sqrt2}\ts f_{12'}.
\een
In the general setting of Sec.~\ref{subsec:mps}, these elements are understood
as the basis $X^1,\dots,X^{14}$.
The elements $X_1,\dots,X_{14}$
of the dual basis with respect to the form \eqref{formg} are then given by the
following expressions, where we use the corresponding capital letters
to think of the $X_i$ as abstract generators of $\U(\g)$ rather than matrices:
\ben
2\tss F_{11}-F_{22}-F_{33},\qquad F_{11}+F_{22}-2\tss F_{33},\qquad
3\tss F_{ji} \quad\text{with}\quad 1\leqslant i,j\leqslant 3\fand i\ne j,
\een
together with
\ben
2\tss F_{41}-\sqrt2\ts F_{23'},\qquad
2\tss F_{42}-\sqrt2\ts F_{31'},\qquad
2\tss F_{43}-\sqrt2\ts F_{12'},
\een
and
\ben
2\tss F_{14}-\sqrt2\ts F_{3'2},\qquad
2\tss F_{24}-\sqrt2\ts F_{1'3},\qquad
2\tss F_{34}-\sqrt2\ts F_{2'1}.
\een
Using \eqref{G}, we can now define
the entries $G_{ij}$ of the matrix $G$ from the expansion
\beql{Ggtwo}
G=\sum_{i,j=1}^7 e_{ji}\ot G_{ij}\in\End\CC^7\ot\U(\g).
\eeq
In particular,
\ben
G_{11}=2\tss F_{11}-F_{22}-F_{33},\qquad
G_{22}=2\tss F_{22}-F_{11}-F_{33},\qquad
G_{33}=2\tss F_{33}-F_{11}-F_{22},
\een
so that $G_{11}+G_{22}+G_{33}=0$. Also, for all
$1\leqslant i,j\leqslant 3$ with $i\ne j$ we have $G_{ij}=3\tss F_{ij}$.
Furthermore,
\ben
G_{14}=2\tss F_{14}-\sqrt2\ts F_{3'2},\qquad
G_{24}=2\tss F_{24}-\sqrt2\ts F_{1'3},\qquad
G_{34}=2\tss F_{34}-\sqrt2\ts F_{2'1},
\een
and
\ben
G_{41}=2\tss F_{41}-\sqrt2\ts F_{23'},\qquad
G_{42}=2\tss F_{42}-\sqrt2\ts F_{31'},\qquad
G_{43}=2\tss F_{43}-\sqrt2\ts F_{12'}.
\een
The remaining entries of the matrix $G$ are determined by the symmetry
properties $G^t=-G$ which give $G_{ij}+G_{j'i'}=0$ together with
\ben
G_{14}=-\sqrt2\ts G_{3'2},\qquad
G_{24}=-\sqrt2\ts G_{1'3},\qquad
G_{34}=-\sqrt2\ts G_{2'1},
\een
and
\ben
G_{41}=-\sqrt2\ts G_{23'},\qquad
G_{42}=-\sqrt2\ts G_{31'},\qquad
G_{43}=-\sqrt2\ts G_{12'}.
\een
Note that the above formulas define an explicit embedding of $\g$ into
the orthogonal Lie algebra $\oa_7$ spanned by the elements $F_{ij}=E_{ij}-E_{j'i'}$,
where the $E_{ij}$ denote the standard basis elements of $\gl_7$.
An expression for the element $\Om$ defined in \eqref{Om}
can be given by
\begin{align}
\Om&=3\ts\sum_{i,j=1}^3\tss f_{ij}\ot f_{ji}-\sum_{i,j=1}^3\tss f_{ii}\ot f_{jj}
+2\ts\sum_{i=1}^3 \tss \big(f_{4i}\ot f_{i4}+ f_{i4}\ot f_{4i}\big)
\non\\[0.3em]
&+\circlearrowleft^{}_{1,2,3}\big(f_{12'}\ot f_{2'1}+f_{2'1}\ot f_{12'}\big)
\label{Omgtwo}\\[0.5em]
&-\sqrt2\ts\circlearrowleft^{}_{1,2,3}\big(f_{14}\ot f_{23'}+ f_{23'}\ot f_{14}
+f_{41}\ot f_{3'2}+f_{3'2}\ot f_{41}\big),
\non
\end{align}
where the symbol $\circlearrowleft^{}_{1,2,3}$ indicates the summation over cyclic
permutations of the indices $1,2,3$ keeping all other symbols, including primes,
at their positions; that is,
\ben
\circlearrowleft^{}_{1,2,3}\ts X_{1,2',3}=X_{1,2',3}+X_{3,1',2}+X_{2,3',1}.
\een

Proposition~\ref{prop:matrrel} provides a matrix form \eqref{matrelg}
of the defining relations of $\U(\g)$ with the elements $G$ and $\Om$
defined in \eqref{Ggtwo} and \eqref{Omgtwo}.

\subsection{A formula for $\Om$ as an element of the centralizer algebra}

By its definition \eqref{Om}, the element $\Om$ can be viewed as an operator
\ben
\Om:V\ot V\to V\ot V.
\een
It is easily seen that this operator commutes with the action of the Lie algebra $\g$
on $V\ot V$ given by
\ben
X\mapsto \pi(X)\ot 1+1\ot \pi(X).
\een
This implies that $\Om$ must be a linear combination of the projections
of $V\ot V$ onto its irreducible components.

As a representation of the Lie algebra $\sll_7$, the tensor product $\CC^7\ot \CC^7$
of two copies of $V=\CC^7$ splits into the direct sum of two irreducible components
\ben
\CC^7\ot \CC^7=\La^2(\CC^7)\oplus \Sr^2(\CC^7),
\een
afforded by the exterior and symmetric square of $V$. The canonical projections
onto the irreducible components are given by the respective operators
$(1-P)/2$ and $(1+P)/2$,
where $P$ is the permutation operator
\ben
P=\sum_{i,j=1}^7 e_{ij}\ot e_{ji}.
\een
The restriction of the representation $\La^2(\CC^7)$ to the subalgebra $\oa_7$
remains irreducible, whereas the restriction of $\Sr^2(\CC^7)$ splits
into two irreducible components; each of them remains irreducible under
the further restriction to the subalgebra $\g\subset\oa_7$ of type $G_2$:
\ben
\Sr^2(\CC^7)\cong V_0\oplus V_{2\om_1}.
\een
The respective projections are given by the operators $Q/7$ and
$(1+P)/2-Q/7$,
where $Q$ is the operator
\ben
Q= \sum_{i,j=1}^7 e_{ij}\ot e_{i'j'}.
\een
It is obtained by applying the antidiagonal transposition
\ben
t:\End\CC^7\to \End\CC^7,\qquad (e_{ij})^t=e_{j'i'},
\een
to the first or the second component of
the permutation operator, $Q=P^{t_1}=P^{t_2}$. Here $V_0$ is the trivial one-dimensional
representation of $\g$, and $V_{2\om_1}$ is the $27$-dimensional representation
corresponding to the double of the first fundamental weight $\om_1$.

When restricted to $\g$, the exterior square
$\La^2(\CC^7)$ splits into two irreducible components
\ben
\La^2(\CC^7)\cong V_{\om_1}\oplus V_{\om_2}
\een
of the respective dimensions $7$ and $14$,
associated with the fundamental weights.
The respective projections are given by the operators $T/6$ and
$(1-P)/2-T/6$,
where the operator $T$ is defined by means of a $3$-form on $\CC^7$ as follows.
The $3$-form $\be$ is defined in the canonical basis $e_1,\dots,e_7$ of $\CC^7$ by
\ben
\be=\sum_{i=1}^3 e_i\wedge e_4\wedge e_{i'}+\sqrt2\ts e_1\wedge e_2\wedge e_3
+\sqrt2\ts e_{3'}\wedge e_{2'}\wedge e_{1'}.
\een
One easily verifies that $\be$ is invariant under the action of $\g$,
which provides a $\g$-module embegging $\CC^7\hra \La^2(\CC^7)$ defined as contraction
with $\be$. Explicitly, the operator $T$ can be written as
\ben
T=\sum_{i,j,k,l=1}^7\sum_{a=1}^7 \be_{ika}\ts\be_{jla}\ts e_{ij}\ot e_{kl},
\een
where the coefficients of the form $\be$ are defined
by
\ben
\be=\sum_{i,j,k=1}^7 \be_{ijk}\ts e_i\ot e_j\ot e_k
\een
via the standard embedding $\La^3(\CC^7)\hra (\CC^7)^{\ot 3}$ such that
\ben
e_{i_1}\wedge e_{i_2}\wedge e_{i_3}=\sum_{\si\in\Sym_3}\sgn\si\cdot
e_{i_{\si(1)}}\ot e_{i_{\si(2)}}\ot e_{i_{\si(3)}}.
\een

Thus, we have an expansion
\ben
\Om=a\ts \Pi_0+b\ts \Pi_{\om_1}+c\ts \Pi_{\om_2}+d\ts \Pi_{2\om_1}
\een
into a linear combination of the projections to the respective irreducible components
of the decomposition
\ben
\CC^7\ot\CC^7\cong V_0\oplus V_{\om_1}\oplus V_{\om_2}\oplus V_{2\om_1}.
\een
The coefficients in the expansion are found by applying $\Om$ to
particular vectors to give the formula
\beql{Ompr}
\Om=1+P-2\tss Q-T.
\eeq

Since $\Pi_0$, $\Pi_{\om_1}$, $\Pi_{\om_2}$ and $\Pi_{2\om_1}$
are pairwise orthogonal projections, we derive the
following relations for the operators $P$, $Q$ and $T$. They commute pairwise and
\beql{pqt}
P^2=1,\qquad Q^2=7\tss Q,\qquad T^2=6\tss T,\qquad P\tss Q=Q,\qquad P\tss T=-T,
\qquad Q\tss T=0.
\eeq

By \eqref{Ompr} the matrix form \eqref{matrelg} of the defining relations
implies a uniform expression for the commutators of the generators of $\g$,
\begin{align}
[G_{ij},G_{kl}]=\de_{kj}\ts G_{il}-\de_{il}\ts G_{kj}&
-2\ts\de_{ki'}\ts G_{j'l}+2\ts\de_{j'l}\ts G_{ki'}
\non\\
{}&+\sum_{a,b=1}^7\Big(\be_{iab}\ts\be_{jlb}\ts G_{ka}-\be_{ikb}\ts\be_{jab}\ts G_{al}\Big).
\label{commrelG}
\end{align}

\bre\label{rem:rmatrix}
It is known by Ogievetsky~\cite{o:fs}
that a rational $R$-matrix associated with $\g$ can be given by the formula
\ben
R(u)=1-\frac{P}{u}+\frac{2\tss Q}{u-6}+\frac{T}{u-4}.
\een
Its expansion into a power series in $u^{-1}$ takes the form
\ben
R(u)=1-(\Om-1)\ts u^{-1}+\dots
\een
so that the classical limit of the $RTT$ relations defining the Yangian $\Y(\g)$
reproduces the defining relations \eqref{matrelg} for $\U(\g)$; see \cite{d:qg}.
\ere

\subsection{Isomorphism with the Chevalley presentation}
\label{subsec:ichp}

The simple Lie algebra $\g$ of type $G_2$ is associated with the Cartan matrix
$A=[a_{ij}]$,
\ben
A=\begin{bmatrix}\phantom{-}2&-1\ts\\
                -3&\phantom{-}2\ts
                \end{bmatrix}.
\een
Its Chevalley presentation is defined by generators
$e_i,h_i,f_i$ with $i=1,2$, subject to
the defining relations
\begin{alignat}{2}
[e_i,f_j]&=\de_{ij}h_i,\qquad [h_i,h_j]&&=0,
\non\\
[h_i,e_j]&=a_{ij}\tss e_j,\qquad [h_i,f_j]&&=-a_{ij}\tss f_j,
\non
\end{alignat}
together with (the Serre relations)
\ben
(\ad\ts e_1)^2\ts e_2=0,\qquad (\ad\ts e_2)^4\ts e_1=0,\qquad
(\ad\ts f_1)^2\ts f_2=0,\qquad (\ad\ts f_2)^4\ts f_1=0.
\een
We let $\al$ and $\be$ denote
the simple roots. The set of positive roots is
\ben
\al,\quad\be,\quad\al+\be,\quad\al+2\tss \be,\quad\al+3\tss \be,\quad 2\tss \al+3\tss \be.
\een
For each positive root $\ga$ we let $e_{\ga}$ and $f_{\ga}$ denote the root
vectors associated with $\ga$ and $-\ga$, respectively. In particular,
$e_1=e_{\al}$ and $e_2=e_{\be}$. We have the triangular decomposition
\beql{triang}
\g=\n_-\oplus\h\oplus\n_+,
\eeq
where the subalgebras $\n_+$ and $\n_-$ are spanned by the positive and
negative root vectors, respectively, while $\h$ is the Cartan subalgebra
with the basis elements $h_1$ and $h_2$.

By the symmetry properties of the matrix $G$ described in Sec.~\ref{subsec:liegtwo},
the Lie algebra $\g$ is spanned by the elements $G_{ij}$ with $1\leqslant i,j\leqslant 4$.
Moreover, due to the relations $G_{11}+G_{22}+G_{33}=0$
and $G_{44}=0$ this leaves a $14$-dimensional Lie algebra. Explicit commutators
between these elements can be obtained by using the embedding $\g\hra\oa_7$
defined in Sec.~\ref{subsec:liegtwo} or by applying the commutation relations
\eqref{commrelG}. Assuming that the indices $i,j,k,l$ run over the set $\{1,2,3\}$,
the commutation relations of $\g$ are then completely described as follows:
\ben
\bal{}
[G_{ij},G_{kl}]&=3\ts\de_{kj}\ts G_{il}-3\ts\de_{il}\ts G_{kj},\\[0.4em]
[G_{ij},G_{k4}]&=3\ts\de_{kj}\ts G_{i4}-\de_{ij}\ts G_{k4},\\[0.4em]
[G_{ij},G_{4l}]&=-3\ts\de_{il}\ts G_{4j}+\de_{ij}\ts G_{4l},\\[0.4em]
[G_{i4},G_{4l}]&=2\ts G_{il},
\eal
\een
together with the relations
\ben
[G_{i4},G_{j4}]=2\sqrt2\ts G_{4k},\qquad [G_{4i},G_{4j}]=-2\sqrt2\ts G_{k4};
\een
the latter hold for the triples $(i,j,k)$ of the form $(1,2,3)$, $(2,3,1)$
and $(3,1,2)$.

\bpr\label{prop:isomche}
The mapping
\ben
e_1\mapsto \frac13\ts G_{12},\qquad e_2\mapsto \frac1{\sqrt2}\ts G_{24},
\qquad f_1\mapsto \frac13\ts G_{21},\qquad f_2\mapsto \frac1{\sqrt2}\ts G_{42}
\een
defines an isomorphism between the Chevalley and the matrix presentations
of $\g$. Moreover, under this isomorphism,
\ben
h_1\mapsto \frac13\ts \big(G_{11}-G_{22}\big)\Fand h_2\mapsto G_{22}.
\een
\epr

\bpf
This follows easily from the commutation relations for the $G_{ij}$.
\epf

By definition, the simple root subspaces $\g_{\ga}$ of $\g$
are $\g_{\al}=\langle G_{12}\rangle$ and $\g_{\be}=\langle G_{24}\rangle$, so that
by Proposition~\ref{prop:isomche} the remaining root subspaces
associated with positives roots
are spanned by the following elements:
\ben
\g_{\al+\be}=\langle G_{14}\rangle,
\qquad \g_{\al+2\be}=\langle G_{43}\rangle,\qquad \g_{\al+3\be}=\langle G_{23}\rangle,
\qquad \g_{2\al+3\be}=\langle G_{13}\rangle.
\een
Similarly, $\g_{-\al}=\langle G_{21}\rangle$ and $\g_{-\be}=\langle G_{42}\rangle$
together with
\ben
\g_{-\al-\be}=\langle G_{41}\rangle,
\qquad \g_{-\al-2\be}=\langle G_{34}\rangle,\qquad \g_{-\al-3\be}=\langle G_{32}\rangle,
\qquad \g_{-2\al-3\be}=\langle G_{31}\rangle.
\een
Therefore, in the triangular decomposition \eqref{triang} we have
\beql{npm}
\n_+=\langle G_{12}, G_{24}, G_{14}, G_{43}, G_{23}, G_{13}\rangle
\fand
\n_-=\langle G_{21}, G_{42}, G_{41}, G_{34}, G_{32}, G_{31}\rangle.
\eeq

\section{Invariants of the vacuum module}
\label{sec:ava}
\setcounter{equation}{0}

Suppose that $\g$ is an arbitrary simple Lie algebra equipped
with an invariant symmetric bilinear form $\langle\ts,\rangle$
which is determined uniquely, up to a nonzero factor.
We will need matrix presentations of the affine Kac--Moody algebra $\wh\g$
analogous to those given in Proposition~\ref{prop:matrrel} for $\g$.
The Lie algebra $\wh\g$
is defined as the central
extension \eqref{km},
where $\g[t,t^{-1}]$ is the Lie algebra of Laurent
polynomials in $t$ with coefficients in $\g$.
The commutation relations
have the form
\beql{defrelaff}
\big[X[r],Y[s]\big]=[X,Y][r+s]+r\ts\de_{r,-s}\langle X,Y\rangle\ts K,
\qquad X, Y\in\g,
\eeq
and the element $K$ is central in $\wh\g$.

\subsection{Matrix presentations of affine Kac--Moody algebras}

We will keep the data associated with the Lie algebra $\g$,
chosen in the beginning of Sec.~\ref{subsec:mps}.
In particular, we will use
a faithful representation \eqref{pi} of $\g$ and the element $\Om$
defined in \eqref{Om}. In addition to the matrix $G$ defined in \eqref{G},
for any $r\in\ZZ$ introduce the matrix $G[r]$ by
\beql{Gr}
G[r]=\sum_{i=1}^d \pi(X^i)\ot X_i[r]\in \End V\ot\U(\wh\g).
\eeq
By analogy with \eqref{Got} introduce the elements
of the algebra $\End V\ot \End V\ot \U(\wh\g)$ by
\ben
G_1[r]=\sum_{i=1}^d \pi(X^i)\ot 1\ot X_i[r]\Fand
G_2[r]=\sum_{i=1}^d 1\ot \pi(X^i)\ot X_i[r].
\een

\bpr\label{prop:matrrelaff}
The defining relations of $\U(\wh\g)$ are equivalent to the matrix relation
\beql{matrelgaff}
G_1[r]\tss G_2[s]-G_2[s]\tss G_1[r]=-\Om\ts  G_2[r+s]+G_2[r+s]\ts \Om
+r\ts\de_{r,-s}\ts\Om\ts K.
\eeq
\epr

\bpf
The proof is completed as for Proposition~\ref{prop:matrrel}, where we use
\eqref{defrelaff} and the additional observation that the sum
\ben
\sum_{i,j=1}^d \langle X_i,X_j\rangle\ts \pi(X^i)\ot \pi(X^j).
\een
coincides with $\Om$.
\epf

\subsection{Segal--Sugawara vectors}
\label{subsec:ssv}

We will now use the notation $\g$ for the
simple Lie algebra of type $G_2$ and follow \cite{f:lc} to recall some general facts
on the affine vertex algebra associated with $\g$.
We equip $\g$ with
the bilinear form \eqref{formg}.

The universal enveloping algebra at the critical level $\U(\wh\g)_{\text{cri}}$
is the quotient of $\U(\wh\g)$
by the ideal generated by $K+12$. Define
the vacuum module $V(\g)_{\text{cri}}$ at the critical level over $\wh\g$
as the quotient of $\U(\wh\g)_{\text{cri}}$ by the left ideal
generated by $\g[t]$. The Feigin--Frenkel center
$\z(\wh\g)$ is defined as the subspace of invariants of the vacuum module; see \eqref{ffvac}.
Any element of $\z(\wh\g)$ is called a Segal--Sugawara vector.
As we pointed out in the Introduction, $\z(\wh\g)$ is a commutative associative algebra
which can be regarded as a subalgebra of $\U\big(t^{-1}\g[t^{-1}]\big)$.
The algebra $\U\big(t^{-1}\g[t^{-1}]\big)$ is equipped with the derivation
$D=-d/dt$ whose action on the generators is given by
\beql{transldef}
D:X[r]\mapsto -r\tss X[r-1],\qquad X\in\g, \quad r<0.
\eeq
The subalgebra $\z(\wh\g)$ of $\U\big(t^{-1}\g[t^{-1}]\big)$
is $D$-invariant. By \cite{ff:ak} (see also \cite{f:lc}) we have the following.

\begin{ffthm}\label{ff:gtwo}
There exist Segal--Sugawara vectors $S_2$ and $S_6$ such that all elements
$D^k S_2$ and $D^k S_6$ with $k\geqslant 0$ are algebraically independent,
and every element of $\z(\wh\g)$ is a polynomial in the $D^k S_2$ and $D^k S_6$.
\qed
\end{ffthm}

We call a pair $S_2$, $S_6$ satisfying the conditions of the theorem
a {\em complete set} of Segal--Sugawara vectors for $\g$.
It is known that the elements $S_2$ and $S_6$ must have the degrees $2$ and $6$
with respect to the canonical filtration of $\U\big(t^{-1}\g[t^{-1}]\big)$
(thus justifying the notation).
These are the degrees of algebraically independent
generators of the algebra of $\g$-invariants of the symmetric algebra $\Sr(\g)$.
The vector $S_2$ of degree $2$ must be proportional to
the canonical Segal--Sugawara vector \eqref{canss},
whereas no formula has been known for a vector of degree $6$. Our first main result
is an explicit formula for such a vector (see Theorem A below).
We keep the notation $G_{ij}$ for the elements of the Lie algebra $\g$ as defined
in Sec.~\ref{sec:liea}. In particular, in accordance with \eqref{Gr},
for any $r\in\ZZ$ we have
\beql{Ggtwoaff}
G[r]=\sum_{i,j=1}^7 e_{ji}\ot G_{ij}[r]\in\End\CC^7\ot\U(\wh\g).
\eeq
Set
\begin{align}\label{sone}
S_2=\tr\ts G[-1]^2&,\\[0.6em]
S_6=\tr\ts G[-1]^6&+5\ts \tr\ts G[-2]\tss G[-1]^4 +14\ts \tr\ts G[-3]\tss G[-1]^3
+456\ts \tr\ts G[-3]^2
\non\\[0.6em]
{}&-639\ts \tr\ts G[-2]^3 +31\ts \tr\ts G[-2]^2\tss G[-1]^2
-312\ts \tr\ts G[-3]\tss G[-2]\tss G[-1].
\label{stwo}
\end{align}

\begin{mthma}
The elements $S_2$ and $S_6$ form a complete set of Segal--Sugawara vectors for
the Lie algebra $\g$.
\end{mthma}

\bpf
Since $\g$ is a simple Lie algebra,
the Lie algebra $\g[t]$ is generated by all elements $G_{ij}[0]$ together with
one nonzero element $G_{kl}[1]$. Therefore, in order to verify the property
$S_2,S_6\in\z(\wh\g)$, it will be sufficient to demonstrate that
in the vacuum module, $G_{ij}[0]\ts S_a=0$ for all $i,j$ and $G_{11}[1]\ts S_a=0$
for $a=2,6$. The first part of these relations will be implied by
the following lemma.

\ble\label{lem:gotr}
For any negative integers $r_1,\dots,r_p$, we have the relations
\ben
G_1[0]\ts\tr\ts G_2[r_1]\dots G_2[r_p]=0
\een
in the algebra
\beql{tenprse}
\End \CC^7\ot \End \CC^7\ot \U(\wh\g),
\eeq
where the elements in $\U(\wh\g)$ are considered
modulo the left ideal $\U(\wh\g)\ts\g[t]$ and
the trace is taken
with respect to the second copy of $\End \CC^7$.
\ele

\bpf
By \eqref{matrelgaff},
\ben
G_1[0]\tss G_2[s]-G_2[s]\tss G_1[0]=-\Om\ts  G_2[s]+G_2[s]\ts \Om.
\een
The argument is completed in the same way as for Corollary~\ref{cor:casim}.
\epf

\ble\label{lem:gone}
For all $k\geqslant 1$ we have the relations
\ben
G_1[1]\ts \tr\ts G_2[-1]^k=\sum_{i=1}^k\tr\ts \Big(G_2[-1]^{i-1}\ts \Om\ts G_2[-1]^{k-i}\ts\Om
-\Om\ts G_2[-1]^{i-1}\ts \Om\ts G_2[-1]^{k-i}\Big)
\een
in the algebra \eqref{tenprse}
modulo the left ideal of $\U(\wh\g)$ generated by $K+12$ and $\g[t]$, where
the trace is taken
with respect to the second copy of $\End \CC^7$.
\ele

\bpf
Applying \eqref{matrelgequi} and \eqref{matrelgaff} we find
\ben
G_1[1]\ts \tr\ts G_2[-1]^k=\sum_{i=1}^k\tr\ts
G_2[-1]^{i-1} \Big(\Om\ts G_1[0]-G_1[0]\ts \Om+K\ts \Om\Big)
G_2[-1]^{k-i}.
\een
Furthermore,
\ben
G_1[0]\ts G_2[-1]^{k-i}=-\Om\ts G_2[-1]^{k-i}+G_2[-1]^{k-i}\ts \Om,
\een
and by applying the partial transposition $t_2$ to both sides we also get
\ben
G_1[0]\ts \big(G_2[-1]^{k-i}\big)^t=-\Om\ts \big(G_2[-1]^{k-i}\big)^t+\big(G_2[-1]^{k-i}\big)^t\ts \Om
\een
since $\Om^{t_2}=-\Om$ as implied by \eqref{Omgtwo}. Hence,
using the general property
\beql{sklyanin}
\tr\ts AB=\tr\ts A^tB^t,
\eeq
by applying $t_2$
to the factors, we get
\ben
\bal
\tr\ts
&G_2[-1]^{i-1}\ts G_1[0]\tss \Om
\tss G_2[-1]^{k-i}=-\tr\ts
\big(G_2[-1]^{i-1}\big)^t\tss G_1[0]
\tss \big(G_2[-1]^{k-i}\big)^t\ts \Om\\[0.5em]
{}&=\tr\ts
\big(G_2[-1]^{i-1}\big)^t\tss \Om
\ts \big(G_2[-1]^{k-i}\big)^t\tss \Om
-\tr\ts
\big(G_2[-1]^{i-1}\big)^t
\tss \big(G_2[-1]^{k-i}\big)^t\ts \Om^2.
\eal
\een
By applying \eqref{sklyanin} again we can write this as
\ben
\tr\ts
\Om\ts G_2[-1]^{i-1}\tss \Om\ts
G_2[-1]^{k-i}
-\tr\ts
G_2[-1]^{i-1}\tss\big(\Om^2\big)^{t_2}
\ts G_2[-1]^{k-i}.
\een
Bringing the calculations together, we obtain
\ben
\bal
G_1[1]\ts \tr\ts G_2[-1]^k&=\sum_{i=1}^k\tr\ts \Big(G_2[-1]^{i-1}\ts \Om\ts G_2[-1]^{k-i}\ts\Om
-\Om\ts G_2[-1]^{i-1}\ts \Om\ts G_2[-1]^{k-i}\Big)\\
{}&+\sum_{i=1}^k\tr\ts G_2[-1]^{i-1}\ts
\Big(\big(\Om^2\big)^{t_2}-\Om^2+K\ts\Om\ts\Big) G_2[-1]^{k-i}.
\eal
\een
Since $P^{t_2}=Q$ and $Q^{t_2}=P$, the relation $\Om^{t_2}=-\Om$ implies
$T^{t_2}=2-P-Q-T$, and so we derive from
\eqref{pqt} that $\big(\Om^2\big)^{t_2}=\Om^2+12\ts \Om$, thus completing the proof.
\epf

It is immediate from Lemma~\ref{lem:gone} that $S_2$ is a Segal--Sugawara vector.
In principle, the lemma can be useful for checking this property of
$S_6$ as well; however, this leads to rather cumbersome expressions
which we were unable to handle without computer's assistance.
We used the Symbolic Manipulation System FORM; see \cite{v:nf}.
Namely, we verified the relation $G_{11}[1]\tss S_6=0$ in the vacuum module
by employing a program within the FORM which works as follows.
We fix a total ordering $\prec$ on the basis elements $G_{ij}[r]$ of $\g$
with $1\leqslant i,j\leqslant 4$ (excluding $(i,j)=(3,3)$ and $(4,4)$)
with the property that $G_{ij}[r]\prec G_{kl}[s]$ for $r<s$.
The input of the program is the element $G_{11}[1]\tss S_6$ written
explicitly as a linear combination of monomials in the generators $G_{ij}[r]$.
The output is a linear combination of the ordered monomials
which is calculated with the use of the commutation relations \eqref{matrelgaff}
written in terms of the $G_{ij}[r]$. The program confirms that
the last factor of each monomial in this linear combination is of the form
$G_{ij}[0]$ or $G_{ij}[1]$.

To complete the proof of the theorem, note that the symbols
$\overline S_2$ and $\overline S_6$ of the elements $S_2$ and $S_6$
in the associated graded algebra $\Sr\big(t^{-1}\g[t^{-1}]\big)\cong\gr\ts\U\big(t^{-1}\g[t^{-1}]\big)$
are given by
\ben
\overline S_2=\tr\ts \overline G[-1]^2\Fand
\overline S_6=\tr\ts \overline G[-1]^6,
\een
where we use the bar notation to indicate
objects associated with the symmetric algebras.
It is easily seen
(e.g., by taking the Chevalley images)
that the elements $\tr\ts \overline G^2$ and $\tr\ts \overline G^{\tss 6}$ are algebraically independent
generators of the algebra of $\g$-invariants $\Sr(\g)^{\g}$.
Thus, the elements $\overline S_2$ and $\overline S_6$ are the respective images
of $\tr\ts \overline G^2$ and $\tr\ts \overline G^{\tss 6}$ under the embedding
$\Sr(\g)\hra\Sr\big(t^{-1}\g[t^{-1}]\big)$
taking $X\in\g$ to $X[-1]$. Therefore, by a general argument of
\cite[Sect.~3.5.1]{f:lc}, the elements $S_2$ and $S_6$ form a complete set
of Segal--Sugawara vectors.
\epf

\bpr\label{prop:moress}
The elements
\ben
\bal
S_3&=\tr\ts G[-1]^3,\\
S_4&=\tr\ts G[-1]^4+\tr\ts G[-2]\tss G[-1]^2,\\
S_5&=\tr\ts G[-1]^5+\tr\ts G[-2]\tss G[-1]^3,
\eal
\een
are Segal--Sugawara vectors for $\g$.
\epr

\bpf
For the element $S_3$ this follows from Lemma~\ref{lem:gone} with $k=3$
by the application of \eqref{sklyanin}. The claim for $S_4$ can also be verified
with the use of Lemma~\ref{lem:gone} and some additional arguments
which we will omit. In fact, for all three elements the claim is also verified
by the same computer program as for $S_6$; see the proof of Theorem A.
\epf

By Theorem A, each of the vectors $S_3$, $S_4$ and $S_5$ is a polynomial
in the $D^k S_2$ and $D^k S_6$. In particular, a direct argument shows that
$S_3=-3\ts D\ts S_2$. Indeed, first we note the easily verified relation
\ben
(G[-1]^2)^t=G[-1]^2+12\ts G[-2],
\een
so that
by applying \eqref{sklyanin} we find
\ben
S_3=\tr\ts (G[-1]^2)^t\ts G[-1]^t=-S_3-12\ts \tr\ts G[-2]\ts G[-1].
\een
This gives $S_3=-6\ts \tr\ts G[-2]\ts G[-1]$ which equals $-3\ts D\ts \tr\ts G[-1]^2$
since
\ben
\tr\ts G[-2]\ts G[-1]=\tr\ts G[-1]\ts G[-2].
\een

We will be able to write down the remaining polynomials after
calculating the Harish--Chandra images of the Segal--Sugawara vectors $S_a$;
see Corollary~\ref{cor:relss} below.

\section{Affine Harish-Chandra isomorphism}
\label{sec:hch}
\setcounter{equation}{0}

It was proved in \cite{ff:ak} that
for any simple Lie algebra $\g$ the algebra $\z(\wh\g)$ is isomorphic
to the classical $\Wc$-algebra $\Wc({}^L\g)$ associated with the Langlands dual
Lie algebra ${}^L\g$ (corresponding to the transposed of the Cartan matrix of $\g$).
We will follow \cite[Sec.~8.1]{f:lc} to recall these results before
applying them to a particular case of type $G_2$.

\subsection{Feigin--Frenkel center and classical $\Wc$-algebra}
\label{subsec:ffcw}

Suppose that $\g=\n_-\oplus\h\oplus \n_+$ is
a triangular decomposition for
a simple Lie algebra $\g$. Regard $\h$ as a subalgebra of $\wh\g$ via the
embedding taking $H\in\h$ to $H[0]$. The adjoint action of $\h$ on $t^{-1}\g[t^{-1}]$
extends to the universal enveloping algebra, and
we have the homomorphism
for the $\h$-centralizer
\beql{hchaff}
\f:\U\big(t^{-1}\g[t^{-1}]\big)^{\h}\to \U\big(t^{-1}\h[t^{-1}]\big)
\eeq
which is the projection to the first summand in the direct sum decomposition
\ben
\U\big(t^{-1}\g[t^{-1}]\big)^{\h}= \U\big(t^{-1}\h[t^{-1}]\big)\oplus
\Jr,
\een
where $\Jr$ is the intersection of the centralizer with the left ideal
of $\U\big(t^{-1}\g[t^{-1}]\big)$ generated by $t^{-1}\n_-[t^{-1}]$.
The Feigin--Frenkel center $\z(\wh\g)$ is a commutative subalgebra
of the centralizer, and
the restriction
of the homomorphism \eqref{hchaff} to
$\z(\wh\g)$ yields an affine version
of the Harish-Chandra isomorphism:
\beql{hchiaff}
\f:\z(\wh\g)\to \Wc({}^L\g),
\eeq
where
the classical $\Wc$-algebra $\Wc({}^L\g)$ is defined as a subalgebra
of $\U\big(t^{-1}\h[t^{-1}]\big)$ which consists of the elements
annihilated by the screening operators; see \cite[Theorem~8.1.5]{f:lc}.
The algebra $\Wc({}^L\g)$ is known to possess algebraically independent
families of generators. Their explicit form in the classical types goes back
to \cite{a:tf} and \cite{gd:fh} in type $A$ and to \cite{ds:la}
in types $B$, $C$ and $D$. For type $G_2$ such generators were calculated in
\cite{mr:cw}; see also \cite{bm:ss}, where they appear in a different context.
The images of explicit generators of $\z(\wh\g)$
under the isomorphism \eqref{hchiaff} are found in \cite{cm:ho} and \cite{ct:qs} in type $A$
and in \cite{mm:yc} for types $B$, $C$ and $D$; see also \cite{r:nf}
for a direct calculation for the Pfaffian-type vector.
In the next section we calculate the images of the Segal--Sugawara vectors
in type $G_2$ constructed in Sec.~\ref{subsec:ssv}.

\subsection{Harish-Chandra images}
\label{subsec:hchi}

Now we assume $\g$ is of type $G_2$ and use its Chevalley presentation; see Sec.~\ref{subsec:ichp}.
Note that the Cartan
matrix $A$ is not symmetric.
Although the Langlands dual ${}^L\g$ is a Lie algebra isomorphic to $\g$,
this leads to a difference between the forms
of the screening operators associated with $A$ and its transpose.
In accordance with \cite[Sec.~8.1.2]{f:lc},
the classical $\Wc$-algebra $\Wc({}^L\g)$
is a subalgebra of $\U(t^{-1}\h[t^{-1}])$ which is
defined as the intersection of the kernels
of the screening operators $V_1$ and $V_2$.
The operators are given by
\beql{scr}
V_i=\sum_{r=0}^{\infty} V_{i\ts r}\sum_{j=1}^2 a_{ji}\ts\frac{\di}{\di\tss h_j[-r-1]},
\eeq
where $h_j[s]=h_j\tss t^s$
and
the coefficients $V_{i\ts r}$ are found by the
relation
\beql{genfv}
\sum_{r=0}^{\infty} V_{i\ts r}\ts z^r=\exp\Big(\sum_{m=1}^{\infty}
\frac{h_i[-m]}{m}\tss z^m\Big).
\eeq

The subalgebra $\Wc({}^L\g)$ is invariant under the action of the derivation $D=-d/dt$
defined in \eqref{transldef}.
This follows from the easily verified relations
\ben
V_i\tss D=\big(D+h_i[-1]\big)\tss V_i,\qquad i=1,2.
\een

To write down explicit generators of $\Wc({}^L\g)$ introduce elements
$g_1[r]$ and $g_2[r]$ of $\U(t^{-1}\h[t^{-1}])$ by
\ben
g_1[r]=h_1[r]+\frac13\ts h_2[r]\Fand g_2[r]=\frac13\ts h_2[r],\qquad r<0.
\een
Moreover, we will regard the derivation $-d/dt$ as a differential operator
and denote it by $\tau$ to distinguish it from the derivation $D$.
In other words, we consider the algebra which is isomorphic
to the tensor product $\U(t^{-1}\h[t^{-1}])\ot\CC[\tau]$ as a vector space, with the relations
\ben
\tau\ts H[r]-H[r]\ts\tau=-r\ts H[r-1],\qquad H\in\h.
\een

A version of the Miura transformation of type $G_2$ was produced in \cite{mr:cw};
cf. \cite{bm:ss}. In the above notation write the product
\begin{multline}
\big(\tau-2g_1[-1]-g_2[-1]\big)\big(\tau-g_1[-1]-2g_2[-1]\big)\big(\tau-g_1[-1]+g_2[-1]\big)\\[0.7em]
{}\times\tau\ts
\big(\tau+g_1[-1]-g_2[-1]\big)\big(\tau+g_1[-1]+2g_2[-1]\big)
\big(\tau+2g_1[-1]+g_2[-1]\big)
\label{miura}
\end{multline}
as a polynomial in $\tau$, so that it equals
\beql{miuraexp}
\tau^7+w_2\tss\tau^5+\dots+w_6\tss\tau+w_7,\qquad w_i\in \U(t^{-1}\h[t^{-1}]).
\eeq
By the results of \cite{mr:cw}, all coefficients $w_i$ belong to $\Wc({}^L\g)$,
and the elements $D^k\tss w_2$ and $D^k\tss w_6$ with $k\geqslant 0$
are algebraically independent generators of $\Wc({}^L\g)$.

The first claim can be verified directly by re-writing the screening operators in terms
of the elements $g_i[r]$. We have the formulas (up to overall scalar factors):
\ben
V_1=\sum_{r=0}^{\infty} V_{1\ts r}\ts\Big(\frac{\di}{\di\tss g_1[-r-1]}-\frac{\di}{\di\tss g_2[-r-1]}\Big),
\een
where
\ben
\sum_{r=0}^{\infty} V_{1\ts r}\ts z^r=\exp\Big(\sum_{m=1}^{\infty}
\frac{g_1[-m]-g_2[-m]}{m}\tss z^m\Big),
\een
and
\ben
V_2=\sum_{r=0}^{\infty} V_{2\ts r}\ts\Big(\frac{\di}{\di\tss g_1[-r-1]}-\frac{2\ts\di}{\di\tss g_2[-r-1]}\Big),
\een
where
\ben
\sum_{r=0}^{\infty} V_{2\ts r}\ts z^r=\exp\Big(\sum_{m=1}^{\infty}
\frac{3\ts g_2[-m]}{m}\tss z^m\Big).
\een
We show that the polynomial \eqref{miura} is annihilated by $V_1$ and $V_2$ by using
the relations
\ben
V_1\ts \tau=\big(\tau+g_1[-1]-g_2[-1]\big)\ts V_1\Fand
V_2\ts \tau=\big(\tau+3\ts g_2[-1]\big)\ts V_2.
\een
In particular, $V_1w_2=V_2w_2=0$, where the coefficient $w_2$ in the expansion \eqref{miuraexp}
is found by
\beql{wtwo}
w_2=-6\ts \big(g_1[-1]^2 + g_1[-1]g_2[-1] + g_2[-1]^2 - 3 g_1[-2] - 2 g_2[-2]\big).
\eeq
We will use the standard notation $w'_2=D\tss w_2$, $w''_2=D^2\tss w_2$, etc. for its
derivatives.
We have the relations
\begin{align}
2\ts w_3 &= 5\ts w'_2,\non\\
4\ts w_4 &= w_2^2 + 12\ts w''_2,
\label{wrel}\\
4\ts w_5 &= 3\ts w_2w'_2 + 8\ts w'''_2.
\non
\end{align}

Our second main result provides the Harish-Chandra images
of the Segal--Sugawara vectors for $\g$ produced in Sec.~\ref{subsec:ssv}.
We use the isomorphism of Proposition~\ref{prop:isomche} to identify the
Chevalley presentation of $\g$ with its matrix presentation.
The triangular decomposition of $\g$ is defined by the subalgebras
$\n_+$ and $\n_-$ given in \eqref{npm}. For the elements of the Cartan
subalgebra we have
\ben
G_{11}[r]\mapsto 3\tss g_1[r]\Fand G_{22}[r]\mapsto 3\tss g_2[r],\qquad r<0.
\een

\begin{mthmb}
For the images under the isomorphism \eqref{hchiaff} we have
\ben
\bal
S_2&\mapsto -6\ts w_2,\\[0.5em]
S_3&\mapsto 18\ts w'_2,\\[0.5em]
S_4&\mapsto 9\ts w_2^2-36\ts w''_2,\\[0.5em]
S_5&\mapsto -63\ts w_2w'_2+72\ts w'''_2,\\[0.2em]
S_6&\mapsto 162\ts w_6-\frac{33}{2}\ts w_2^3+\frac{63}{2}\ts (w'_2)^2
+90\ts w_2w''_2-576\ts w^{(4)}_2.
\eal
\een
\end{mthmb}

\bpf
The first image is easily verified using the commutations relations
of $\wh\g$. Since
\eqref{hchiaff} is a differential algebra isomorphism, the second image
follows from the formula $S_3=-3\ts S'_2$, where we extend the derivative notation
$S'=D\ts S$ to elements of $\z(\wh\g)$. The remaining relations are verified by
a computer program within the FORM; see \cite{v:nf}. Namely, the Segal--Sugawara vectors
$S_a$ are written as linear combinations of monomials in the generators
$G_{ij}[r]$. The program first provides the respective linear combinations of ordered
monomials, where the total ordering on the generators is chosen in a way
consistent with the triangular decomposition of $\g$. The resulting
Harish-Chandra images $\f(S_a)$ are linear combinations of the monomials containing only
the diagonal generators $G_{11}[r]$ and $G_{22}[r]$.
The program provides an explicit linear combination
for $\f(S_a)$ in terms of monomials in the generators of $\Wc({}^L\g)$.
\epf

\bre
For verification purposes, we used another program to make sure that
all images $\f(S_a)$ with $a=2,\dots,6$ do belong to the classical $\Wc$-algebra $\Wc({}^L\g)$.
This was done by applying the screening operators to check that
$V_1\ts\f(S_a)=V_2\ts\f(S_a)=0$.
\qed
\ere

\bco\label{cor:relss}
We have the relations
\begin{align}
S_3 &= -3\ts S'_2,\non\\[0.3em]
4\ts S_4 &= S_2^2 + 24\ts S''_2,
\non\\[0.3em]
4\ts S_5 &= -7\ts S_2S'_2 - 48\ts S'''_2.
\non
\end{align}
\eco

\bpf
This is immediate from Theorem B and relations \eqref{wrel}.
\epf

\section{Commutative subalgebras and Gaudin model}
\label{sec:agh}
\setcounter{equation}{0}

Here we consider applications of Theorems A and B to
explicit constructions of maximal commutative subalgebras of $\U(\g)$
and to the Gaudin model associated with the simple Lie algebra $\g$ of type $G_2$.

\subsection{Quantization of the shift of argument subalgebras}

We will use the matrix presentation of the Lie algebra $\g$ so that
it is spanned by the entries $G_{ij}$ of the matrix $G$ given in \eqref{Ggtwo}.
For any element $\mu\in\g^*$ we set $\mu_{ij}=\mu(G_{ij})\in\CC$ so that we can regard
$\mu=[\mu_{ij}]$ as a matrix of the form
\ben
\mu=\sum_{i,j=1}^7 e_{ji}\ot \mu_{ij}.
\een
The bilinear form \eqref{formg} allows us to identify $\g^*$ with $\g$.
An element $\mu\in\g^*\cong\g$ is called {\em regular},
if the centralizer $\g^{\mu}$ of $\mu$ in $\g$ has minimal possible
dimension. This minimal dimension coincides with the rank of $\g$ and so equals $2$.

The next theorem is implied by Theorem A and
a positive solution of Vinberg's quantization problem~\cite{v:sc}
given by Rybnikov~\cite{r:si}
and Feigin, Frenkel and Toledano Laredo~\cite{fft:gm}
with the use of the algebra $\z(\wh\g)$; cf. \cite{fm:qs} and \cite{m:ff}
for such applications in classical types.

For a given arbitrary element $\mu$
introduce polynomials in a variable $z$ with coefficients in $\U(\g)$ by the formulas
\ben
\bal
A(z)&=\tr\tss(G+\mu\tss z)^2,\\[0.4em]
B(z)&=\tr\tss(G+\mu\tss z)^6+5\ts\tr\ts G\tss(G+\mu\tss z)^4
+14\ts\tr\ts G\tss(G+\mu\tss z)^3+31\ts\tr\ts G^2\tss(G+\mu\tss z)^2,
\eal
\een
and write
\ben
\bal
A(z)&=A_0\tss z^2+A_1\tss z+A_2,\\[0.3em]
B(z)&=B_0\tss z^6+B_1\tss z^5+B_2\tss z^4+B_3\tss z^3+B_4\tss z^2+B_5\tss z+B_6.
\eal
\een

\begin{mthmc}
For any $\mu\in\g^*$, all elements $A_i$ and $B_j$ of $\U(\g)$ pairwise commute.
Moreover, if $\mu$ is regular, then the elements
$A_1,A_2,B_1,B_2,B_3,B_4,B_5,B_6$
are algebraically
independent and generate a maximal commutative subalgebra of $\U(\g)$.
\end{mthmc}

\bpf
Given
a nonzero $z\in\CC$, the mapping
\beql{evalr}
\varrho^{}_{\ts\mu,z}:\U\big(t^{-1}\g[t^{-1}]\big)\to \U(\g),
\qquad G_{ij}[r]\mapsto G_{ij}\tss z^r+\de_{r,-1}\ts\mu_{ij},
\eeq
defines an algebra homomorphism.
The Feigin--Frenkel center $\z(\wh\g)$ is a commutative
subalgebra of $\U\big(t^{-1}\g[t^{-1}]\big)$
so that the image of $\z(\wh\g)$ under $\varrho^{}_{\ts\mu,z}$
is a commutative subalgebra of $\U(\g)$.
This image is independent of $z$ and
will be denoted by $\Ac_{\mu}$.
For the images of the generator matrices we have
\ben
z\tss\varrho^{}_{\ts\mu,z}\big(G[-1]\big)=G+\mu\ts z
\Fand z^r\varrho^{}_{\ts\mu,z}\big(G[-r]\big)=G\quad\text{for}\quad r\geqslant 2.
\een
Consider the elements $S_2,S_6\in\z(\wh\g)$
provided by Theorem A. All coefficients of the polynomials in $z$ defined by
\ben
z^2\varrho^{}_{\ts\mu,z}(S_2)\Fand
z^6\varrho^{}_{\ts\mu,z}(S_6)
\een
belong to the commutative subalgebra $\Ac_{\mu}$. Moreover, $A(z)=z^2\varrho^{}_{\ts\mu,z}(S_2)$,
while $z^6\varrho^{}_{\ts\mu,z}(S_6)$ equals $B(z)$
plus a linear combination of the terms $\tr\ts G^3$, $\tr\ts G^2$ and
$\tr\ts G^2\tss(G+\mu\tss z)$.
However, using \eqref{sklyanin} and the easily verified relation $(G^2)^t=G^2+12\ts G$,
we get
\ben
\tr\ts G^3=-6\ts \tr\ts G^2\Fand \tr\ts G^2\tss(G+\mu\tss z)=-6\ts \tr\ts G\tss(G+\mu\tss z).
\een
Therefore, the constant term of the linear combination is proportional to $A_2$, whereas
the coefficient of $z$ is proportional to $A_1$.
This proves the first part of the theorem.

The second part will follow by considering the symbols $\overline A_i$
and $\overline B_j$ of the coefficients
$A_i$ and $B_j$ in the associated graded algebra $\gr\ts\U(\g)\cong \Sr(\g)$.
They are given by
\ben
\bal
\tr\tss(\overline G+\mu\tss z)^2&=\overline A_0\tss z^2+\overline A_1\tss z+\overline A_2,\\[0.3em]
\tr\tss(\overline G+\mu\tss z)^6&=\overline B_0\tss z^6+\overline B_1\tss z^5
+\overline B_2\tss z^4+\overline B_3\tss z^3
+\overline B_4\tss z^2+\overline B_5\tss z+\overline B_6.
\eal
\een
It was shown in \cite{fft:gm} that
if the element $\mu$ is regular, then the coefficients
$\overline A_1,\overline A_2,\overline B_1,\dots,\overline B_6$
are algebraically independent generators of the Poisson commutative
subalgebra $\overline\Ac_{\mu}$ of $\Sr(\g)$, known
as the {\it Mishchenko--Fomenko subalgebra\/} or
{\it shift of argument subalgebra\/}.
This also follows from the earlier results of Bolsinov~\cite{b:cf}.
Furthermore, it was shown in
\cite{py:as} that $\overline\Ac_{\mu}$ is maximal Poisson commutative.
This implies the second statement of the theorem.
\epf

\bre\label{rem:computer}
As an additional test, we employed a computer program within the FORM~\cite{v:nf}
(similar to the one used in the proof of Theorem~A)
to verify some particular cases of Theorem~C. Namely, we considered
diagonal matrices
$
\mu=\diag\ts [\mu_1,\mu_2,\mu_3,0,-\mu_3,-\mu_2,-\mu_1]
$
with $\mu_1+\mu_2+\mu_3=0$ and regarded $\mu_1$ and $\mu_2$ as variables.
The elements $A_2$ and $B_6$ are central, while
$A_1$ and $B_1$ belong to the Cartan subalgebra of $\g$ and
so commute with all other elements $B_i$. Our program verified that $B_2$ commutes
with each of $B_3$, $B_4$ and $B_5$, and that $B_3$ commutes with $B_4$ and $B_5$.
\qed
\ere

\subsection{Eigenvalues of the Gaudin Hamiltonians}

A connection of the center at the critical level $\z(\wh\g)$
with the Gaudin Hamiltonians was first observed by
Feigin, Frenkel and Reshetikhin~\cite{ffr:gm}. They used the
Wakimoto modules over the affine Kac--Moody algebra $\wh\g$ to calculate
the eigenvalues of the Hamiltonians on the Bethe vectors of the Gaudin model
associated with an arbitrary simple Lie algebra $\g$.
Given an element $S\in\z(\wh\g)$, the eigenvalues are expressed in terms
of the Harish-Chandra image of $S$; see also more recent work \cite{fft:gm}
for generalizations to non-homogeneous
Hamiltonians.

We will use the explicit formulas for elements of $\z(\wh\g)$ provided
by Theorem A and the general result of \cite[Theorem~6.7]{fft:gm}
to write explicit Gaudin operators and their
eigenvalues on Bethe vectors for the simple Lie algebra $\g$ of type $G_2$.

Using coassociativity of the standard coproduct on $\U\big(t^{-1}\g[t^{-1}]\big)$
defined by
\ben
\Delta: G_{ij}[r]\mapsto G_{ij}[r]\ot 1+1\ot G_{ij}[r],\qquad r<0,
\een
for any $\ell\geqslant 1$ we get the homomorphism
\beql{comult}
\U\big(t^{-1}\g[t^{-1}]\big)\to \U\big(t^{-1}\g[t^{-1}]\big)^{\ot\tss\ell}
\eeq
as an iterated coproduct map. Now fix distinct complex numbers
$z_1,\dots,z_\ell$ and let $u$ be a complex parameter.
Applying homomorphisms of the form \eqref{evalr} to the tensor factors in
\eqref{comult}, we get another homomorphism
\beql{psiu}
\Psi:\U\big(t^{-1}\g[t^{-1}]\big)\to \U(\g)^{\ot\tss\ell},
\eeq
given by
\ben
\Psi:G_{ij}[r]\mapsto \sum_{a=1}^\ell (G_{ij})_a(z_a-u)^r+\de_{r,-1}\ts\mu_{ij}\in \U(\g)^{\ot\tss\ell},
\een
where $(G_{ij})_a=1^{\ot (a-1)}\ot G_{ij}\ot 1^{\ot (\ell-a)}$; see \cite{r:si}.
We will twist this homomorphism by the involutive anti-automorphism
\beql{sgn}
\vs: \U\big(t^{-1}\g[t^{-1}]\big)\to \U\big(t^{-1}\g[t^{-1}]\big),
\qquad G_{ij}[r]\mapsto -G_{ij}[r],
\eeq
to get the anti-homomorphism
\beql{phiu}
\Phi:\U\big(t^{-1}\g[t^{-1}]\big)\to \U(\g)^{\ot\tss\ell},
\eeq
defined as the composition $\Phi=\Psi\circ\vs$.
Since the Feigin--Frenkel center
$\z(\wh\g)$ is a commutative subalgebra of $\U\big(t^{-1}\g[t^{-1}]\big)$,
the image of $\z(\wh\g)$ under $\Phi$ is a commutative subalgebra
$\Ac(\g)_{\mu}$ of $\U(\g)^{\ot\tss\ell}$. It depends on the chosen
parameters $z_1,\dots,z_\ell$, but does not depend on $u$
\cite{r:si}; see also \cite[Sec.~2]{fft:gm}.

For any $\la\in\h^*$, the Verma module $M_{\la}$ is defined as the quotient of $\U(\g)$
by the left ideal generated by $\n_+$ and the elements $h_i-\la(h_i)$
with $i=1,2$. We denote the image of $1$ in $M_{\la}$ by $1_{\la}$.
Given weights $\la_1,\dots,\la_\ell\in\h^*$ consider the tensor product
of the Verma modules $M_{\la_1}\ot\dots\ot M_{\la_\ell}$. Formulas for common
eigenvectors (the {\em Bethe vectors})
\ben
\phi(w_1^{i_1},\dots,w_m^{i_m})\in M_{\la_1}\ot\dots\ot M_{\la_\ell}
\een
for the commutative subalgebra $\Ac(\g)_{\mu}$
in this tensor product can be found in \cite{ffr:gm} (also
reproduced in \cite{mm:eb}). They are parameterized by
a set of distinct complex numbers $w_1,\dots, w_m$
with $w_i\ne z_j$ and a collection (multiset) of labels $i_1,\dots,i_m\in\{1,2\}$.

Suppose that $\mu\in\h^*$. We regard $\mu$ as a functional on
$\g$ which vanishes on $\n_+$ and $\n_-$. The system
of the {\em Bethe ansatz equations} takes the form
\beql{bae}
\sum_{i=1}^\ell\frac{\la_i(h_{i_j})}{w_j-z_i}-\sum_{s\ne j}
\frac{\al_{i_s}(h_{i_j})}{w_j-w_s}=\mu(h_{i_j}),
\qquad j=1,\dots,m,
\eeq
where $\al_1=\al$ and $\al_2=\be$ are the simple roots; cf. \cite{bm:ss}.

Introduce the homomorphism
from $\U\big(t^{-1}\h[t^{-1}]\big)$ to rational functions in $u$ by the rule:
\beql{hw}
\varrho:G_{ii}[-r-1]\mapsto \frac{\di_u^{\ts r}}{r!}\ts \Gc_{i}(u),\qquad
r\geqslant 0,\qquad i=1,2,
\eeq
where
\ben
\Gc_{i}(u)=\sum_{a=1}^\ell \frac{\la_a(G_{ii})}
{u-z_a}-\sum_{j=1}^m\frac{\al_{i_j}(G_{ii})}{u-w_j}-\mu_{ii}.
\een

Consider the Segal--Sugawara vectors $S_a$ with $a=2,\dots,6$ provided by Theorem~A and
Proposition~\ref{prop:moress}. Their Harish-Chandra images
$\f(S_a)\in \U\big(t^{-1}\h[t^{-1}]\big)$
are given in Theorem~B.
The composition $\varrho\circ\f$ takes each Segal--Sugawara vector $S_a$ to
the rational function $\varrho\big(\f(S_a)\big)$ in $u$.
Furthermore,
we regard the image $\Phi(S)$ of $S$ under the anti-homomorphism \eqref{phiu}
as an operator in the tensor product of Verma modules
$M_{\la_1}\ot\dots\ot M_{\la_\ell}$.
The following is a consequence of \cite[Theorems~6.5 and 6.7]{fft:gm};
cf. \cite{mm:eb}.

\begin{mthmd}
Suppose that the Bethe ansatz equations \eqref{bae} are satisfied.
If the Bethe vector $\phi(w_1^{i_1},\dots,w_m^{i_m})$ is nonzero, then it is an eigenvector
for the operator $\Phi(S_a)$ with the eigenvalue $\varrho\big(\f(S_a)\big)$
for each $a=2,\dots,6$.
\qed
\end{mthmd}

\end{document}